\input amstex
\documentstyle{amsppt}
\magnification 1200
\NoRunningHeads
\NoBlackBoxes
\document

\def\tg_+{\tilde{\frak g_+}}

\def\ad{\text{ad}}
\def\Hom{\text{Hom}}

\def\a{\frak a}

\def\Ua{U_q(\tilde\g)}
\def\U2{{\Ua}_2}
\def\g{\frak g}

\def\Z{\Bbb Z}

\def\Q{\Bbb Q}

\def\l{\lambda}
\def\I{\Cal I}

\def\<{\langle}
\def\>{\rangle}
\def\o{\otimes}
\def\e{\varepsilon}

\def\End{\text{End}}

\def\b{\frak b}

\topmatter
\title Quantization of Lie bialgebras, II
\endtitle
\author {\rm {\bf Pavel Etingof and David Kazhdan} \linebreak
\vskip .1in
Department of Mathematics\linebreak
Harvard University\linebreak 
Cambridge, MA 02138, USA\linebreak
e-mail: etingof\@math.harvard.edu\linebreak kazhdan\@math.harvard.edu}
\endauthor
\endtopmatter

\centerline{\bf Part III} 
\centerline{\bf Classification of quantized universal enveloping algebras}
\vskip .1in
\centerline{\bf Abstract}
\vskip .03in
This paper is a continuation of \cite{EK}.
We show that the quantization procedure of \cite{EK} is given by
universal acyclic formulas and defines a functor from the category
of Lie bialgebras to the category of quantized universal enveloping algebras.
We also show that this functor defines an equivalence
between the category of Lie bialgebras over $k[[h]]$
and the category quantized universal enveloping (QUE) algebras.
\vskip .1in

\centerline{\bf Acknowledgements}

The work of the first author was partically supported 
by an NSF postdoctoral fellowship. 

\vskip .1in

\centerline{\bf 1. Universality and functoriality of the quantization of
Lie bialgebras.}
\vskip .1in

In this paper we will use the notation of \cite{EK}.

\vskip .1in 

1.1. {\it Linear algebraic structures.}

Here we introduce the notion of a linear algebraic structure 
which is borrowed from \cite{La}.

 Let $(\Cal C,\o)$ be a symmetric
monoidal $\Q$-linear category whose objects are nonnegative
integers, such that $[n]=[1]^{\o n}$ (the unit object is $[0]$). 
We will call such a category a cyclic category 
(by analogy with a cyclic group).

Let $\Cal C$ be a cyclic category.
Let $S=\cup_{m,n\ge 0} S_{mn}$ be a set of morphisms of
$\Cal C$, $S_{mn}\subset\Hom_{\Cal C}([m],[n])$, $m,n\in \Z_+$. 
We say that $\Cal C$ is generated by $S$ 
if any morphism of $\Cal C$ can be obtained from
the morphisms in $S_{mn}$ 
and permutation morphisms in $\Hom_{\Cal C}(m,m)$ by iterating
three elementary 
operations: 

1) composition of morphisms;

2) tensor product of morphisms;

3) linear combinations of morphisms over $\Q$.

Let $Z=\cup Z_{mn}$, $m,n\in\Z_+$, be a bigraded set.

\proclaim{Definition} The free cyclic category
 $\Cal F_Z$ is a cyclic category equipped with
a collection of maps $\mu_{mn}: Z_{mn}\to \Hom_{\Cal F_Z}(m,n)$
such that

(i) $\Cal F_Z$ is generated by $\cup_{m,n\ge 0}\mu_{mn}(Z_{mn})$, and

(ii) for any cyclic category $\Cal C$ equipped with a collection of maps
$\phi_{mn}: Z_{mn}\to \Hom_{\Cal C}([m],[n])$ there exists 
a unique symmetric tensor functor (ACU functor, see \cite{DM})
$F:\Cal F_Z\to \Cal C$ such that $F([1]_{\Cal F_Z})=[1]_{\Cal C}$, and
$F(\mu_{mn}(z))=\phi_{mn}(z)$, $z\in Z_{mn}$. 
\endproclaim   

The pair $(\Cal F_Z,\{\mu_{mn}\})$ with such properties 
exists and is unique up to symmetric tensor equivalence.

{\bf Remark.} Roughly, the free cyclic category is a cyclic category
generated by a set of morphisms without any nontrivial relations.

Let $(\Cal C,\o)$ be a $\Q$-linear monoidal category, and $\I=\{\I_{X,Y}
\subset \Hom_{\Cal C}(X,Y),X,Y\in
\Cal C\}$ be a collection of subspaces. We say that a morphism
$\phi\in \Hom_{\Cal C}(X,Y)$ is $\I$-negligible if $\phi\in \I_{X,Y}$.
We say that $\I$ is a tensor ideal in 
$\Cal C$ if the composition in any order
of any morphism with an $\I$-negligible morphism
 is $\I$-negligible, and the tensor product of any morphism
with an $\I$-negligible morphism is $\I$-negligible. 

Given a tensor ideal $\I$ in 
$\Cal C$. Define the quotient category $\Cal D=\Cal C/\I$ as follows. 
The objects of $\Cal D$ are the same as those of $\Cal C$, and
$\Hom_{\Cal D}(X,Y)=\Hom_{\Cal C}(X,Y)/\I_{X,Y}$. 
It is easy to show that $\Cal D$
carries a natural structure of a $\Q$-linear monoidal category. We will
say that an identity between morphisms holds in $\Cal C$ modulo $\I$
if it holds in $\Cal D$. 

Let $J=\{J_{X,Y}
\subset \Hom_{\Cal C}(X,Y),X,Y\in
\Cal C\}$, be a collection of subsets.
We denote by $\<J\>$ the smallest tensor ideal in $\Cal C$ such 
that $J_{X,Y}\subset \<J\>_{X,Y}$ 
for any objects $X,Y$, and say that
 $\<J\>$ is the tensor ideal generated by $J$. 

\proclaim{Proposition 1.1} Let $\Cal C$ be any cyclic category
generated by a set $S$ of morphisms. Then $\Cal C$ has
the form $\Cal F_S/\I$, where $\I$ is a tensor ideal in $\Cal F_S$. 
\endproclaim

\demo{Proof} Easy.$\square$\enddemo

Let $\Cal C$ be a cyclic category, $\I$ be a tensor ideal in $\Cal C$.
For any positive integer $n$,
define the power $\I^n$ to be the tensor ideal in $\Cal C$ generated
by the elements $\phi_1\circ ...\circ \phi_n$, where $\phi_1,...,\phi_n\in \I$.

Now we define the notion of completion of a cyclic category 
with respect to a tensor ideal. Let $\Cal C$ be a cyclic category,
$\I$ be a tensor ideal in $\Cal C$. Then the spaces
of morphisms of the cyclic categories
$\Cal C/\I^n$ form a projective system. Let
$\Cal C_\I$ be the cyclic category whose objects are the same as
those of $\Cal C$, and
$\Hom_{\Cal C/\I}(X,Y)=\underleftarrow{\lim} \Hom_{\Cal C/\I^n}(X,Y)$, 
$X,Y\in\Cal C$. This category is 
called the completion of $\Cal C$ with respect to $\I$. 

Define a topological cyclic category to be a cyclic category
in which the sets of morpisms are topological vector spaces over $\Q$,
and composition of morphisms is continuous. The category
$\Cal C_\I$ has a natural structure of a topological cyclic category, 
where the topology is defined by the ideal $\I$.

Throughout this Chapter,
let $\Cal N$ be a symmetric monoidal $\Q$-linear category, and
$X$ be an object in $\Cal N$.

Let $\Cal C$ be a cyclic category.

\proclaim{Definition} 
A linear algebraic structure of type $\Cal C$ on $X$ is 
a symmetric tensor functor $G:\Cal C\to\Cal N$ such that $G([1])=X$. 
\endproclaim

Thus, a linear algebraic structure of type $\Cal C$ on $X$ 
is a collection of morphisms between tensor
powers of $X$ which satisfy certain consistency relations.
If an object $X$ equipped with such morphisms has been fixed, the functor $G$
corresponding to it is denoted by $G_X$. 

The same definition applies to the case when $\Cal C$ and $\Cal N$ are 
topological categories.

Denote by $\Cal G(\Cal C,\Cal N)$
the category of symmetric tensor fucntors from $\Cal C$ to $\Cal N$,
i.e. the category of linear algebraic structures of type $\Cal C$ 
on objects of $\Cal N$.

Let $\Cal C$ be a (topological) cyclic category. 
We say that $\Cal C$ is nondegenerate if the natural map
$\Bbb Q[S_n]\to \Hom([n],[n])$ is injective. From now on, we consider
only nondegenerate categories.
 
Let $\hat\Cal C$ be the category obtained by formal addition 
to $\Cal C$ of the kernels of all the idempotents $P\in \Q[S_n]$ 
acting on the objects $[n]$, $n\ge 1$.  
The objects of $\hat C$ 
are pairs $([n],P)$, where $n\ge 0$, and $P\in \Bbb Q[S_n]\to
\Hom_{\Cal C}([n],[n])$
is an idempotent, where morphisms are defined by
$\Hom_{\hat\Cal C}(([n],P),([m],Q))=\{f\in\Hom_{\Cal C}([n],[m]): 
f\circ P=Q\circ f=f\}$. Let $\bar\Cal C$ be the
closure of $\Cal C$ under inductive limits.
In this category, every object is isomorphic to a direct sum of 
indecomposable ones, and indecomposable objects
correspond to irreducible representations
of $S_n$, $n\ge 1$. In particular, we have objects
$S^n[1]=([n],\text{Sym}_n)$, where $\text{Sym}_n$ is the symmetrizer
in $\Bbb Q[S_n]$.

Let $\Cal N$ be closed under inductive limits. 
Then any linear algebraic structure $G$
of type $\Cal C$ on $X$ extends to an additive symmetric tensor functor
$G: \bar\Cal C\to \Cal N$.

\vskip .05in
1.2. {\it Examples of linear algebraic structures.} 

Some common examples of linear algebraic structures are:

A. Associative algebras with unit.
In this case the set $S$ consists of an element of bidegree (2,1)
(``the universal product''), and an element of bidegree (0,1)
(``the unit''). 
The category $\Cal C=\text{AA}$
is $\Cal F_S/\I$, where $\I$ is generated by the 
the associativity identity for the product and the unit axiom.
An associative algebra  in $\Cal N$ is an object with a linear
algebraic structure of type $\text{AA}$.

B. Lie algebras.
In this case the set $S$ consists of one element of bidegree (2,1)
(``the universal commutator''), 
and the category $\Cal C=\text{LA}$
is $\Cal F_S/\I$, where $\I$ is generated by two relations -- 
skew-symmetry and the Jacobi identity for the commutator.
A Lie algebra  in $\Cal N$ is an object with a linear
algebraic structure of type $\text{LA}$.

We will deal with the
following examples of linear algebraic structures.

1. Lie bialgebras.
In this case the set $S$ consists of two elements of bidegrees (2,1) and (1,2)
(``the universal commutator and cocommutator''), 
and the category $\Cal C=\text{LBA}$
is $\Cal F_S/\I$, where $\I$ is generated by five relations -- 
skew-symmetry and the Jacobi identity for the commutator and cocommutator
and the condition that cocommutator is a 1-cocycle.
A Lie bialgebra  in $\Cal N$ is an object with a linear
algebraic structure of type $\text{LBA}$.

2. Quasitriangular Lie bialgebras.
In this case the set $S$ consists of two elements of bidegrees (2,1) and (0,2)
(``the universal commutator and classical r-matrix''), 
and the category $\Cal C=\text{QTLBA}$
is $\Cal F_S/\I$, where $\I$ is generated by four relations -- 
skew-symmetry and the Jacobi identity for the commutator,
invariance of $r+r^{op}$, and the classical Yang-Baxter equation. 
A quasitriangular Lie bialgebra  in $\Cal N$ is an object with a linear
algebraic structure of type $\text{QTLBA}$.

3. Hopf algebras.
In this case the set $S$ consists of six
 elements of bidegrees (2,1),(1,2),(0,1),(1,0),(1,1),(1,1)
(``the universal product, coproduct, unit, counit, antipode,inverse 
antipode''), 
and the category $\Cal C=\text{HA}$
is $\Cal F_S/\I$, where $\I$ is generated by the relations coming from
the axioms of a Hopf algebra. 
A Hopf algebra in $\Cal N$ is an object with a linear
algebraic structure of type $\text{HA}$.

4. Quasitriangular Hopf algebras. 
In this case the set $S$ consists of eight
 elements of bidegrees (2,1),(1,2),(0,1),(1,0),(1,1),(1,1),(0,2),(0,2)
(``the universal product, coproduct, unit, counit, antipode, inverse antipode,
R-matrix,
inverse R-matrix''), 
and the category $\Cal C=\text{QTHA}$
is $\Cal F_S/\I$, where $\I$ is generated by the relations coming from
the axioms of a quaitriangular Hopf algebra. 
A quasitriangular
Hopf algebra in $\Cal N$ is an object with a linear
algebraic structure of type $\text{QTHA}$.

5. Classical Yang-Baxter algebras. 
In this case the set $S$ consists of three
 elements of bidegrees (2,1),(0,1),(0,2)
(``the universal product, unit, and r-matrix''), 
and the category $\Cal C=\text{CYBA}$
is $\Cal F_S/\I$, where $\I$ is generated by the associativity
relation and the classical Yang-Baxter equation. 
A classical
Yang-Baxter algebra in $\Cal N$ is an object
 with a linear
algebraic structure of type $\text{CYBA}$.

6. Quantum Yang-Baxter algebras. 
In this case the set $S$ consists of three
 elements of bidegrees (2,1),(0,1),(0,2)
(``the universal product, unit, and R-matrix''), 
and the category $\Cal C=\text{QYBA}$
is $\Cal F_S/\I$, where $\I$ is generated by the associativity
relation and the quantum Yang-Baxter equation. 
A quantum 
Yang-Baxter algebra in $\Cal N$ is an object
 with a linear
algebraic structure of type $\text{QYBA}$.

7. Co-Poisson Hopf algebras \cite{Dr1}.
In this case the set $S$ consists of six
 elements of bidegrees (2,1),(1,2),(0,1),(1,0),(1,1),(1,2)
(``the universal product, coproduct, unit, counit, antipode, Poisson 
cobracket''), 
and the category $\Cal C=\text{CPHA}$
is $\Cal F_S/\I$, where $\I$ is generated by the relations coming from
the axioms of a co-Poisson Hopf algebra. 
A co-Poisson Hopf algebra in $\Cal N$ is an object with a linear
algebraic structure of type $\text{CPHA}$.

8. Quasitriangular co-Poisson Hopf algebras. 
In this case the set $S$ consists of six
 elements of bidegrees (2,1),(1,2),(0,1),(1,0),(1,1),(0,2)
(``the universal product, coproduct, unit, counit, antipode, r-matrix''), 
and the category $\Cal C=\text{QTCPHA}$
is $\Cal F_S/\I$, where $\I$ is generated by the relations coming from
the axioms of a quaitriangular co-Poisson Hopf algebra. 
A quasitriangular co-Poisson 
Hopf algebra in $\Cal N$ is an object with a linear
algebraic structure of type $\text{QTCPHA}$.

\vskip .05in

1.3. {\it The functor of universal quantization.}

Let $\Cal C_1$, $\Cal C_2$ be cyclic categories.

\proclaim{Definition} 
A universal construction is a symmetric tensor functor
$Q: \Cal C_2\to\bar\Cal C_1$.
\endproclaim
  
The functor $Q$ defines a functor $\widehat{Q}:\Cal G(\Cal C_1,\Cal N)\to 
\Cal G(\Cal C_2,\bar \Cal N)$.
 
{\bf Examples.} 1. Let $\Cal C_1=AA$, $\Cal C_2=LA$,
and $Q: LA\to AA$, such that $Q([1])=[1]$,
$Q([,])=*-*^{op}$, where $*$ is the product and $*^{op}$ is the opposite 
product. This functor gives rise to a functor $\widehat{Q}$ from the category of
associative algebras in $\Cal N$ to the category of Lie algebras in $\Cal N$,
which assigns to every associative algebra itself
regarded as a Lie algebra.

2. Let $\Cal C_1=LA$, $\Cal C_2=AA$, and $Q: AA\to LA$ be the 
functor such that $\widehat{Q}$ assigns to every Lie algebra $\g$ in $\Cal N$
its universal 
enveloping algebra $U(\g)$ in $\tilde\Cal N$. 
The functor $Q$ satisfies $Q([1])=S[1]$, and
$Q(*)$ is described in terms of the commutator using 
the standard method of computing the product of monomials. 
For example, $Q(*)|_{[1]\o [1]}=\text{Sym}_2+[,]/2$, where
$[1]\o [1]$ is the subobject in $S[1]\o S[1]$ obtained by tensoring
two subobjects $[1]\subset S[1]$.
 
3. The functor $Q$ of taking the universal enveloping algebra
(example 2) extends to a functor $Q: HA\to LA$, since
the universal enveloping algebra of a Lie algebra carries
a natural structure of a Hopf algebra. Moreover, it 
extends to a functor $Q: CPHA\to LBA$, so that the corresponding
functor $\widehat{Q}$ assigns to any Lie bialgebra its universal enveloping 
algebra regarded as a co-Poisson Hopf algebra. 

The main result of Chapter 1 is the following theorem, 
which establishes the universality of the quantization
of Lie bialgebras, $r$-matrices, quasitriangular and triangular 
Lie bialgebras, obtained in the previous chapters. 

Let $(\a,[,],\delta)$ be a Lie bialgebra over $k$.
Define a Lie bialgebra $\a_h$ over $k[[h]]$ to be
$(\a[[h]],[,],h\delta)$. Similarly, if $(\a,[,],r)$ is a quasitriangular
Lie bialgebra over $k$, we define $\a_h=(\a[[h]],[,],hr)$, and if
$(A,*,r)$ is a classical Yang-Baxter algebra, we define $A_h=(A[[h]],*,hr)$. 
 
\proclaim{Theorem 1.2} 

 There exist ``universal quantization functors'' 

(i) 
$Q: \text{HA}_{\<\Delta-\Delta^{op},S-S^{-1}\>}\to 
\overline{\text{LBA}_{\<\delta\>}}$ such that 
for any Lie bialgebra $\a$ over $k$ $\widehat{Q}(\a_h)=U_h(\a)$;

(ii) 
$Q^{qt}: \text{QTHA}_{\<\Delta-\Delta^{op},R-1,S-S^{-1}\>}\to 
\overline{\text{QTLBA}_{\<r\>}}$ such that 
for any quasitriangular 
Lie bialgebra $\a$ over $k$ $\widehat {Q^{qt}}(\a_h)=U_h^{qt}(\a)$,
where $U_h^{qt}(\a)$ 
is the quasitriangular quantization defined in Section 6.1;

(iii) 
$Q^{YB}: \text{QYBA}_{\<R-1\>}\to 
\text{CYBA}_{\<r\>}$ such that 
for any classical Yang-Baxter algebra $(A,r)$ 
over $k$ one has $\widehat {Q^{YB}}(A_h)=(A,R)$, where
$R$ is constructed from $r$ as explained in Chapter 5.
\endproclaim

In the language of Drinfeld \cite{Dr3},
Theorem 1.2 implies that the multiplication and comultiplication in $U_h(\a)$ 
are expressed via the commutator and cocommutator in $\a$ in terms of acyclic 
tensor calculus, and the Hopf algebra relations in $U_h(\a)$ can be 
formally deduced from the axioms of a Lie bialgebra. Thus it answers 
positively Question 1.2 in \cite{Dr3}
(the existence of universal quantization of Lie bialgebras), 
and hence question 2.1 (the existence of
a quantum Campbell-Hausdorff series). Together with the material of 
Chapters 5,6, it also answers positively the ``universality'' questions in 
Sections 3,4 of \cite{Dr3}.
  
The functoriality of quantization implies the existence of quantization
functors over local Artinian or pro-Artinian algebras. Namely,
let $K$ be a commutative local Artinian or pro-Artinian $\Q$-algebra, $I$ be 
the maximal ideal in $K$, $k=K/I$ be the residue field
(main example: $K=k[[h]]$). 
Let $\Cal A_K$ be the category of topologically 
free $K$-modules, i.e. modules of the form $\underleftarrow{\lim}V\o K/I^n$,
where $V$ is a free abelian group. 

Let:

 $LBA_0(K)$ be the category of Lie bialgebras in $\Cal A_K$ 
with $\delta=0$ mod $I$;

 $HA_0(K)$ be the category of QUE algebras, i.e. 
Hopf algebras in $\Cal A_K$ which are cocommutative mod $I$;

 $QTLBA_0(K)$ the category of quasitriangular Lie bialgebras in $\Cal A_K$ 
with $r=0$ mod $I$;

 $QTHA_0(K)$ the category of
quasitriangular QUE algebras, i.e. quasitriangular
Hopf algebras in $\Cal A_K$ which are cocommutative and have $R=1$ 
mod $I$;

 $CYBA_0(K)$ be
the category of classical Yang-Baxter algebras in $\Cal A_K$ 
with $r=0$ mod $I$;

 $QYBA_0(K)$ be 
the category of
quantum Yang-Baxter algebras in $\Cal A_K$ with $R=1$ 
mod $I$.

Then:

(i) The functor $Q$ defines a functor $\widehat{Q}$ 
from $LBA_0(K)$ to $HA_0(K)$;

(ii) The functor $Q^{qt}$ defines a functor $\widehat {Q^{qt}}$ 
from $QTLBA_0(K)$ to $QTHA_0(K)$;

(ii) The functor $Q^{YB}$ defines a functor $\widehat {Q^{YB}}$ 
from $CYBA_0(K)$ to $QYBA_0(K)$. 

{\bf Remark.} Quantization of triangular Lie bialgebras is also 
universal and functorial. Formulations and proofs of the results are the 
same as in the quasitriangular case.

It is imprortant to remember that the functors $Q,Q^{qt},Q^{YB},
\widehat{Q},\widehat {Q^{qt}},\widehat {Q^{YB}}$ depend on the choice 
of the Lie associator $\Phi$, which was fixed in Section 1.3.
We will later show that different choices of $\Phi$ give isomorphic 
functors, but this is a non-trivial fact. Therefore, when we want to
emphasize the dependence of these functors on $\Phi$, we will 
write them as $Q_\Phi,Q^{qt}_\Phi,Q^{YB}_\Phi,
\widehat {Q_\Phi},\widehat {Q_\Phi^{qt}},\widehat {Q_\Phi^{YB}}$.  

The proof of Theorem 1.2 is given in the next section. 
The idea of the proof is to show that that the construction of quantization
described in \cite{EK} actually defines
some universal formulas which can be used to construct the quantization
in a general tensor category. 

In section 1.5, we give an application of functoriality. 
\vskip .1in 

1.4. {\it Proof of Theorem 1.2}

Part (i). 

Let $\a_+$ be
the canonical Lie bialgebra $[1]$ in the tensor category
$\Cal C=\overline{LBA_{\<\delta\>}}$,
 with commutator $\mu$ and cocommutator 
$\delta$.
Let $U(\a_+):=S\a_+\in\Cal C$ be the universal 
enveloping algebra of $\a_+$. Our goal 
is to introduce a Hopf algebra structure
on $U(\a_+)$, which coincides with the standard one modulo $\<\delta\>$,
and yields the Lie bialgebra structure on $\a_+$ when considered 
modulo $\<\delta\>^2$. This Hopf algebra will be exactly $Q([1])$,
where $[1]$ is the generating object of the category $HA_{\<\Delta-\Delta^{op},
S-S^{-1}\>}$.

We will perform virtually the same constructions as in Part II 
of \cite{EK}, remembering, though, that now $\a_+$ is not a vector space
and even not a set but an object in the category.

First, we define the notion of a module, comodule, and
dimodule over a Lie algebra, coalgebra, and bialgebra, respectively,
in a general tensor category. The notion of a dimodule over a Lie bialgebra
is equivalent, for a finite dimensional Lie bialgebra over a field,
 to the notion of a module over its double.

\proclaim{Definition} 

(i) Let $\a_+$ be a Lie algebra in a tensor category 
$\Cal N$ with commutator $\mu$. An object $X\in\Cal N$ is said to be equipped
with the structure of a left $\a_+$-module if it is endowed
with a morphism $\pi: \a_+\o X\to X$ (the action of $\a_+$ on $X$), 
such that $\pi \circ (1\o \pi)=\pi\circ (\mu\o 1)$ on $\Lambda^2\a_+\o X$. 

(ii) Let $\a_+$ be a Lie coalgebra in a tensor category 
$\Cal N$ with cocommutator $\delta$. 
An object $X\in\Cal N$ is said to be equipped
with the structure of a right $\a_+$-comodule if it is endowed
with a morphism $\pi^*: X\to \a_+\o X$ (the coaction of $\a_+$ on $X$), 
such that $\text{Alt}_{21}\circ
(1\o \pi^*)\circ \pi^*=(\delta\o 1)\circ \pi^*$, where 
$\text{Alt}_{21}$ is the alternator of the first and second components
($\text{Alt}(a\o b):=b\o a-a\o b$). 

(iii) Let $\a_+$ be a Lie bialgebra in a tensor category 
$\Cal N$. An object $X\in\Cal N$ is said to be equipped
with the structure of a $\a_+$-dimodule if it is endowed
with two morphisms $\pi: \a_+\o X\to X$, $\pi^*: X\to\a_+\o X$,
such that $\pi$ is a left action of $\a_+$ on $X$ as a Lie algebra,
$\pi^*$ is a right coaction of $\a_+$ on $X$ as a Lie coalgebra,
and they agree according to the formula
$$
\pi^*\circ \pi=(\pi\o 1)\circ \sigma_{12}\circ (1\o \pi^*) -
(1\o \pi)\circ (\delta\o 1)+(\mu\o 1)\circ (1\o \pi^*),\tag 1.1
$$      
where $\sigma_{ij}$ denotes the permutation of the components $i,j$.
\endproclaim

Let $X$ be any $\a_+$-module (comodule, dimodule). Define $X_0$ to be the
object $X$ with the zero structure of a module (comodule, dimodule).

There is an obvious notion of tensor product of modules and comodules.
The tensor product of dimodules
is just the tensor product of the underlying modules and comodules. 
Thus, modules, comodules,
and dimodules over $\a_+$ in $\Cal N$ form a tensor category. 
We denote the first category by $\Cal M_{\a_+}$, the second by
$\Cal M^{\a_+}$, and the third by $\Cal M_{\a_+}^{\a_+}$.

Now we define the Verma dimodules $M_-$, $M_+^*$ over $\a_+$. 
As objects of $\Cal C$, $M_-=S\a_+$, $M_+^*=\hat S\a_+$.
Let $1_-: \bold 1\to M_-$, $1_+^*:\bold 1\to M_+^*$ be 
embeddings of the identity subobjects into $M_-, M_+^*$, and
$1_-^*: M_-\to\bold 1$, 
$1_+: M_+^*\to \bold 1$ be the corresponding projections. 

The action of $\a_+$ in $M_-$ is the same as the standard left action
of $\a_+$ in $U(\a_+)$. The coaction of $\a_+$ in $M_-$ is then 
completely determined (via formula (1.1))
by the coaction on the identity component $\bold 1$ of $M_-$, which we
define to be zero. 

The formulas for the action
(respectively, coaction) of $\a_+$ in $M_+^*$ 
are obtained from the formulas for 
the coaction (respectively, action) in $M_-$ by interchanging $\mu$ and 
$\delta$, reversing the order of compositions, and changing the signs. 

Let $M_1,M_2$ be $\a_+$-dimodules.
Define the classical r-matrix and the Casimir operator in 
$\End_{\Cal C}(M_1\o M_2)$ by 
$$
r=(\pi\o 1)\circ \sigma_{12} \circ (1\o\pi^*), \Omega=r+r^{op}.\tag 1.2
$$

This implies that if $\Phi$ is any Lie associator, then for any three 
$\a_+$-dimodules $M_1,M_2,M_3$ we can define an invertible element
of $\End_{\Cal C}(M_1\o M_2\o M_3)[[h]]$, which represents the action
of $\Phi$ in $M_1\o M_2\o M_3$.

Now we describe 
an intertwiner $\psi$, which is used in the construction of quantization. 

For any object $X\in\Cal M_{\a_+}^{\a_+}$, define the map 
$\theta: \Hom_{\Cal M_{\a_+}^{\a_+}}( M_-\o X_0,M_+^*\o X)\to 
\Hom_{\Cal C}(X_0,X)$ by $\theta(f)=(1_+\o 1_X)\circ f\circ 
(1_- \o 1_{X_0})$. 

\proclaim{Lemma 1.3} The map $\theta$ is an isomorphism.
\endproclaim

\demo{Proof} 
By Frobenius reciprocity, for any two dimodules $X,Y$ over 
$\a_+$,
$\Hom_{\Cal M_{\a_+}^{\a_+}}(M_-\o X,M_+^*\o Y)=
\Hom_{\Cal M^{\a_+}}(X,M_+^*\o Y)=
\Hom_{\Cal C}(X,Y)$. This implies the Lemma. 
\enddemo

Now set $X=M_-$. 
Denote by $\psi$ the morphism $\theta^{-1}(1_X)$. 
Let $\eta=\psi\circ (1_- \o 1_{X_0}): X_0\to M_+^*\o X$. 

Consider the category $\Cal C_h$, with the same objects as in $\Cal C$,
and morphisms defined by $\Hom_{\Cal C_h}(X,Y)=\Hom_{\Cal C}(X,Y)[[h]]$. 
We define the Hopf algebra $U_h(\a_+)$ in $\Cal C_h$
as follows. By the definition,
this is the object $M_-$, with the product $m: M_-\o M_-\to M_-$
 and coproduct $\Delta: M_-\to M_-\o M_-$
defined as follows.

The product is
$$
m^{op}=(1_+\o 1_+\o 1)\circ (1\o\psi) \circ \Phi\circ (\eta\o 1).
\tag 1.3
$$

The coproduct is
$$
\Delta=J^{-1}\circ \Delta_0,\tag 1.4
$$
where $\Delta_0$ is the standard coproduct in $M_-$, and
$J:M_-\o M_-\to M_-\o M_-$ is given by
$$ 
 J=(1_+\o 1\o 1_+\o 1)\circ \Phi_{1,3,24}^{-1}\circ \Phi_{3,2,4} 
\circ e^{-h\Omega_{23}/2}
\circ \Phi_{2,3,4}^{-1}\circ \Phi_{1,2,34}\circ (\eta\o\eta).\tag 1.5
$$
(cf. formula (3.1)).

\proclaim{Proposition 1.4} Formulas (1.3)-(1.5) define a Hopf algebra 
structure on $M_-$. 
\endproclaim

\demo{Proof} The axioms of a Hopf algebra are demonstrated similarly
to the case of Lie bialgebras in the category of vector spaces, 
which was done in \cite{EK}. For example, 
associativity of multiplication is shown by a computation in Chapter 9;
associativity of comultiplication is shown by the computation given
in the proof of Proposition 8.2. 
\enddemo

Define the functor $P_h: \Cal C\to \Cal C_h$
between topological categories, which maps 
objects to themselves, and $P_h(\mu)=\mu, P_h(\delta)=h\delta$.
The image of this functor is the closed 
subcategory $\Cal C_h'\subset \Cal C_h$,
generated over $\Q$ by the morphisms $\mu,h\delta$, and the functor 
$P_h: \Cal C\to \Cal C_h'$ is an equivalence. 

\proclaim{Proposition 1.5} 
The morphisms $m,\Delta$ defined by (1.3)-(1.5) belong to $\Cal C_h'$. 
\endproclaim

\demo{Proof} It is clear that a morphism $f$ of $\Cal C_h$ belongs to 
$\Cal C_h'$ if and only if it is invariant under the transformation
$a^t$ of $\Cal C_h$ defined by 
$h\to th$, $\mu\to\mu$, $\delta\to \delta/t$, for all $t\in \Q^*$.
Therefore, it is enough to show that 
the Hopf algebra structure on $M_-$ defined by (1.3)-(1.5) 
is invariant under the automorphism $a^t$.

Denote by $\a_+^t$ the Lie bialgebra $a^t(\a_+)$.
Let $\Cal M^t$ be the category of $\a_+^t$-dimodules in $\Cal C_h$. 

We have an equivalence of categories $\theta_t^*: \Cal M^1\to\Cal M^t$, 
which maps any $\a_+$-dimodule $X$ to $X$, with an action and coaction of
$\a_+^t$ defined by the formulas
$\pi_t=\pi,\pi_t^*=t\pi^*$. 
In the case of Lie bialgebras over a field, this 
equivalence comes from the natural Lie algebra isomorphism between 
the doubles $D(\g_+^t)$, $D(\g_+)$,
$\theta_t: D(\g_+^t)\to D(\g_+)$, which is the identity on $\g_+$ and 
multiplication by $t$ on $\g_-$. 
By the definition, this functor maps
the Verma modules to the Verma modules, 
sends the classical r-matrix
$r$ to $t^{-1}r$, and preserves the morphism $\psi$.

Since the associator $\Phi$ is a function of $h\Omega_{12},h\Omega_{23}$, 
it is invariant under the transformation $a^t$. By Proposition 1.4, this 
implies that $U_h(\a_+)$ is fixed by $a^t$, and Proposition 1.5 follows. 
\enddemo

Let
$m',\Delta'$ be the morphisms in the category $\Cal C$ such
that $P_h(m')=m$, $P_h(\Delta')=\Delta$.
They exist by Proposition 1.5 and are obviously unique.  
The morphisms $m',\Delta'$ define a Hopf algebra structure on $M_-$ in
$\Cal C$. 

Now define $Q([1])$ to be the object $S[1]$ with the Hopf algebra structure 
$(m',\Delta')$. This determines the functor $Q$ whose existence
is claimed in Part (i) of Theorem 1.2. Part (i) is proved.

Now let us prove parts (ii) and (iii) of Theorem 1.2.
 
Let $\Cal N$ be a symmetric tensor category, $X\subset \Cal N$ an object.
By an element of $X$ we mean a morphism $x: \bold 1\to X$
(we will write $x\in X$).  

Part (ii). Let $\g=[1]$ be the canonical quasitriangular Lie bialgebra in the
 tensor category
$\Cal C=\overline{QTLBA_{\<r\>}}$. 
Let $U(\g)$ be its universal enveloping algebra. We have an element
$r\in \g\o\g$ -- the classical r-matrix of $\g$. 
Our purpose is to define a Hopf algebra $H$ in $\Cal N$, which
coincides with $U(\g)$ as an object of $\Cal N$.

The definition of $H$ is as follows. The product and unit in $H$ are the same 
as in $U(\g)$. The coproduct is defined by the formula
$$
\Delta=J^{-1}\Delta_0J,\tag 1.6
$$
where $J=J(r)$ is an invertible element of $U(\g)\o U(\g)$ (cf. (3.1),(3.2)). 

The element $J(r)$ is defined by a universal formula which is obtained from 
(3.1): $J(r)=1+r/2+O(\<r\>^2)$. 

Define the functor $Q^{qt}$ by $Q^{qt}([1])=H$.
This functor satisfies the conditions of Part (ii) of Theorem 1.2. 
Part (ii) is proved. 

Part (iii). Let $A$ be the canonical 
classical Yang-Baxter algebra in the tensor category $CYBA_{<r>}$, 
and $r\in A\o A$ be the classical r-matrix. Let $J(r)$ be the
element of $A\o A$ given by the same formula as in Part (ii). 
Define $R=(J^{op})^{-1}e^{(r+r^{op})/2}J\in A\o A$.
(cf. (3.10)). As shown in Part I of \cite{EK}, this 
element satisfies the quantum Yang-Baxter equation.   

Define the functor $Q^{YB}$ by $Q^{YB}([1])=(A,R)$.
This functor satisfies the conditions of Part (iii) of Theorem 1.2. 
Part (iii) is proved. 
 
This completes the proof of Theorem 1.2.
\vskip .1in

1.5. 
{\it Identification of two quantizations of a quasitriangular Lie bialgebra.}

Let $\a$ be a finite dimensional
quasitriangular Lie bialgebra over $k$. Let $U_h(\a)$ be the quantization 
of $\a$ constructed in Chapter 4, and $U_h^{qt}(\a)$ be the quasitriangular 
quantization of $\a$ constructed in Section 6.1.

\proclaim{Theorem 1.6} The quantized universal enveloping algebras $U_h(\a)$, 
$U_h^{qt}(\a)$ are isomorphic.
\endproclaim

To prove Theorem 1.6, we first need the following result, which appears
(in somewhat different form) in \cite{RS}.

\proclaim{Lemma 1.7} Let $\a$ be a quasitriangular Lie bialgebra, and $\g$ be the double 
of $\a$. Then the linear map $\tau:\g\to\a$ defined by
$$
\tau(x+f)=x+(f\o 1)(r), x\in\a, f\in\a^*,\tag 1.6
$$
is a homomorphism of quasitriangular Lie bialgebras.
\endproclaim

\demo{Proof} First we show that $\tau$ is a homomorphism of Lie algebras,
i.e. $\tau([g_1g_2])=[\tau(g_1)\tau(g_2)]$. 
This is obvious when $g_1,g_2\in\a$.
Assume that $f,g\in\a^*$. Then, using the classical Yang-Baxter equation, we get
$$
\gather
\tau([fg])=([fg]\o 1)(r)=(f\o g\o 1)((\delta\o 1)(r))=\\
(f\o g\o 1)([r_{13}+r_{23},r_{12}])=(f\o g\o 1)([r_{13},r_{23}])=[\tau(f)\tau(g)].
\tag 1.7\endgather
$$  
Now assume that $x\in\a,f\in \a^*$. Then
$$
\gather
\tau([xf])=\tau(\ad^*x(f))-\tau(\ad^*f(x))=\\
\tau((f\o 1)([r,x\o 1]))+
\tau((f\o 1)([x\o 1+1\o x,r]))=\\
\tau((f\o 1)([1\o x,r]))=[\tau(x)\tau(f)].\tag 1.8
\endgather
$$

Now we check that $\tau$ is a homomorphism of 
quasitriangular Lie bialgebras. Let $\tilde r$ be the quasitriangular structure on $\g$.
If $x_i$ is a basis of $\a$, and $f_i$ is the dual basis of $\a^*$, then 
$\tilde r$ is given by the formula $\tilde r=\sum_i x_i\o f_i$. 
Thus we have
$$
(\tau\o\tau)(\tilde r)=\sum_i\tau(x_i)\o\tau(f_i)=\sum_i x_i\o (f_i\o 1)(r)=r.\tag 1.9
$$
The Lemma is proved.
$\square$\enddemo

\demo{Proof of Theorem 1.6} Lemma 1.8 claims that there exists a 
morphism of quasitriangular Lie bialgebras $\tau:\g\to\a$ which is 
the identity on $\a$. Theorem 1.2 states that 
quasitriangular quantization of Section 6.1 is a functor from the category
of quasitriangular Lie bialgebras to the category of quasitriangular
topological Hopf algebras over $k[[h]]$. Thus, $\tau$ defines a morphism
$\hat\tau: U_h^{qt}(\g)\to U_h^{qt}(\a)$. On the other hand, $U_h(\a)$ was constructed
as a subalgebra in $U_h^{qt}(\g)$, so we have 
an embedding $\eta: U_h(\a)\to U_h^{qt}(\g)$. Consider the morphism
$\tau\circ \eta: U_h(\a)\to U_h^{qt}(\a)$. This morphism is an isomorphism since it equals
to $1$ modulo $h$. The theorem is proved.
$\square$\enddemo

\proclaim{Corollary 1.8} The quantization of the double $\g$ of a 
finite dimensional Lie bialgebra $\a$ constructed in Chapter 3 
is isomorphic to the quantization of $\g$ as a Lie bialgebra, constructed in 
Chapter 4.
\endproclaim
  
{\bf Remark.} The analog of Theorem 1.6 holds for infinite dimensional 
Lie bialgebras. 
Namely, the ``usual'' quantization of $\a$ defined in Section 9 is isomorphic 
to its quasitriangular quantization. The proof is analogous to the finite 
dimensional case.  

\vskip .05in

\centerline{\bf 2. Dequantization of QUE algebras} 
\vskip .03in
 
2.1. {\it The main result.}

The main result of this paper is the following theorem. 

\proclaim{Theorem 2.1} The functor $\widehat{Q}$ is an equivalence 
of categories.
\endproclaim

{\bf Remark.} We plan to prove 
in a forthcoming paper that the functors $\widehat{Q^{qt}}$,
$\widehat{Q^{YB}}$ are also equivalences of categories. 

We will prove the theorem for $K=k[[h]]$. It is not 
difficult to generalize the proof to the case when $K$ is a general
pro-Artinian algebra over $k$. 

In order to prove the theorem, we construct the functors of 
dequantization (quasiclassical limit)
which are inverse to the functors of quantization
$\widehat{Q}$,
$\widehat {Q^{qt}}$, $\widehat {Q^{YB}}$.
The usual notion
of quasiclassical limit described in Chapter 3 is not sufficient for us
since it assigns to a QUE algebra over $k[[h]]$
 a Lie bialgebra over $k$, not over $k[[h]]$, and thus 
erases nearly all the information about the QUE algebra. 
Our construction of quasiclassical limit is different
and uses the action of the Grothendieck-Teichmuller semigroup 
on braided structures on a tensor category, defined by Drinfeld.
This construction is described in the next sections. 
\vskip .05in

2.2. {\it The Grothendieck-Teichmuller semigroup.}

In this section we follow \cite{Dr4}.

Let $\Cal B$ be a $\Q$-linear topological braided tensor category, 
with a tensor ideal $\Cal I$ such that $\Cal B$ is complete 
with respect to the topology defined by $\Cal I$. This means,
any series of morphisms of the form $\sum_{n=0}^\infty a_n$, where 
$a_n\in \Cal I^n$, is convergent (for example, $\Cal B$ is $K$-linear,
and $\Cal I_{X,Y}=I\Hom_{\Cal B}(X,Y)$,
 where $I$ is the maximal ideal in $K$). 
Let $\Phi$ be the associativity
morphism, and $\beta$ the commutativity morphism (braiding) of $\Cal B$.  

\proclaim{Definition} We say that $\Cal B$ is quasisymmetric if
$\beta^2=1\text{ mod } \Cal I$. 
\endproclaim

Denote by $\overline{GT}$ the set of all pairs $a=(\l,f)$ 
such that 

1) $\l\in \Q$,  $f(X,Y)$ is a formal
series of the form $\exp(P(\ln X,\ln Y))$,
where $P$ is a Lie formal series over $\Q$, 
and

2) for any quasisymmetric category with associativity morphism
$\Phi$ and braiding $\beta$ the morphisms
$$
a(\Phi)=\Phi f(\beta_{21}\beta_{12},\Phi^{-1}\beta_{32}\beta_{23}\Phi),\
a(\beta)=\beta\circ(\beta^2)^m, m=\frac{1}{2}(\l-1)\tag 2.1
$$
define a new structure of a quasisymmetric category on $\Cal B$.

The set $\overline{GT}$ is a monoid under
the composition law
$$
\gather
(\l_1,f_1)(\l_2,f_2)=(\l_1\l_2,f_1\circ_{\l_2} f_2),\\
f_1\circ_{\l} f_2(X,Y):=f_1(f_2(X,Y)X^{\l}f_2(X,Y)^{-1},Y^{\l})f_2(X,Y).\tag 2.2\endgather
$$
The subset $GT=\{(\l,f)\in \overline{GT}:\l\ne 0\}$ is a
group under the same composition law. 

For $\mu\in \Q$, denote by $GT_\mu$ the set of elements
of $\overline{GT}$ of the form $(\mu,f)$. 
It is clear that $\overline{GT}$, $GT$, $GT_\mu$
have a natural structure of proalgebraic $\Q$-varieties.  

The following result is due to Drinfeld (\cite{Dr4}, p. 855). 

\proclaim{Theorem 2.2} For any $a_0\in GT_0$ there exists  
a unique algebraic homomorphism of semigroups $a: \Q\to \overline{GT}$, 
such that $a(0)=a_0$, and $a(\mu)\in GT_\mu$, $\mu\in \Q$.
\endproclaim
\vskip .1in

{\bf Remark.} An algebraic homomorphism means a homomorphism, 
which is polynomial modulo any finite degree. 

Given an element $a=(\l,f)\in \overline{GT}$ and a quasisymmetric
tensor category $\Cal B$ with associativity isomorphism $\Phi$ and
braiding $\beta$, we can define a new quasisymmetric category 
$\Cal B^a$, which is the same as $\Cal B$ as an additive category
with tensor product, and the associativity isomorphism and braiding given by
formula (2.1). This category is symmetric if $a\in GT_0$. 
Also, if $a,b\in\overline{GT}$, then $\Cal B^{ab}=(\Cal B^b)^a$.
Thus, the semigroup $\overline{GT}$ acts on the set of 
structures of a quasisymmetric category on $\Cal B$. 
\vskip .1in

2.3. {\it Dimodules over a Hopf algebra.}

In this section we will introduce the notion of a dimodule
over a Hopf algebra $H$ in an arbitrary symmetric tensor category.
It is analogous to the notion of a dimodule over a Lie bialgebra,
introduced in Chapter 1, and, for the case of finite-dimensional Hopf 
algebras, is
equivalent to the notion of a module over the quantum double $D(H)$.

\proclaim{Definition} 

(i) Let $A_+$ be an associative algebra in a symmetric tensor category 
$\Cal N$, with product $m$ and unit $\iota$. 
An object $X\in\Cal N$ is said to be equipped
with the structure of a left $A_+$-module if it is endowed
with a morphism $\pi: A_+\o X\to X$ (the action of $A_+$ on $X$), 
such that $\pi \circ (1\o\pi)=\pi\circ (m\o 1)$ on $A_+\o A_+\o X$,
and $\pi\circ (\iota\o 1)=id$ on $X$.  

(ii) Let $A_+$ be a coassociative coalgebra in a 
symmetric tensor category 
$\Cal N$, with coproduct $\Delta$ and counit $\e$. 
An object $X\in\Cal N$ is said to be equipped
with the structure of a right $A_+$-comodule if it is endowed
with a morphism $\pi^*: X\to A_+\o X$ (the coaction of $A_+$ on $X$), 
such that $(1\o\pi^*)\circ \pi^*= (\Delta^{op}\o 1)\circ \pi^*$ on $X$,
and $(\e\o 1)\circ \pi^*=id$ on $X$. 

(iii) Let $A_+$ be a Hopf algebra in a 
symmetric tensor category 
$\Cal N$. An object $X\in\Cal N$ is said to be equipped
with the structure of an $A_+$-dimodule if it is endowed
with two morphisms $\pi: A_+\o X\to X$, $\pi^*: X\to A_+\o X$,
such that $\pi$ is a left action of $A_+$ on $X$ as an algebra,
$\pi^*$ is a right coaction of $A_+$ on $X$ as a coalgebra,
and they agree according to the formula
(cf \cite{Dr1}, p. 816)
$$
\pi^*\circ \pi=(m_3 \o\pi)
\circ \sigma_{13}\sigma_{24}\circ (S^{-1}\o 1^{\o 4})
\circ (\Delta_3\o\pi^*),\tag 2.3 
$$      
where $m_3:=m\circ (m\o 1)$,
and $\Delta_3:=(\Delta\o 1)\circ\Delta$.
\endproclaim

Let $A_+$ be a Hopf algebra in $\Cal N$. 
We say that an $A_+$-module (comodule) $X$ is trivial if
$\pi=\e\o 1_X$ (respectively, $\pi^*=\iota\o 1_X$). An $A_+$-dimodule
is called trivial if it is trivial both as a module and a comodule. 
For any $A_+$-module, comodule, or dimodule $X$, let $X_0$ be the object $X$
equipped with the trivial structure of a module (comodule, dimodule). 

There is an obvious notion of tensor product of modules and comodules.
Namely, for any two modules (comodules) $V,W$
$$
\pi_{V\o W}=(\pi_V\o\pi_W)\circ\sigma_{23}\circ (\Delta\o 1\o 1);
\pi_{V\o W}^*=(m^{op}\o 1\o 1)\circ \sigma_{23}\circ (\pi_V^*\o\pi_W^*).
\tag 2.4
$$
The tensor product of dimodules
is just the tensor product of the underlying modules and comodules. 
It follows from \cite{Dr1}, p. 816,
and can be checked by a direct computation, that in this way one indeed
obtains a new dimodule.  

Thus, modules, comodules,
and dimodules over $A_+$ in $\Cal N$ form a tensor category. 
We denote the first category by $\Cal M_{A_+}$, the second by
$\Cal M^{A_+}$, and the third by $\Cal M_{A_+}^{A_+}$.

According to the results of Drinfeld \cite{Dr1}, the category 
$\Cal M_{A_+}^{A_+}$ has a natural structure of a braided tensor category.
The braiding is defined by the formula
$$
\beta=\sigma\circ R, R=(\pi\o 1)\circ \sigma_{12}\circ (1\o\pi^*).\tag 2.5
$$
Drinfeld proved that (2.5) satisfies the hexagon relations. 

{\bf Remark.} The existence of this braiding corresponds to the fact
that the double of a Hopf algebra is a quasitriangular Hopf algebra.   

Now we define the Verma dimodules $M_-$, $\hat M_+^*$ over $A_+$. 
As objects of $\Cal N$, $M_-=\hat M_+^*=A_+$. Let $1_-:\bold 1\to M_-$,
$1_+^*: \bold 1\to \hat M_+^*$ be the unit of $A_+$, and 
$1_+: \hat M_+^*\to\bold 1$, $1_-^*: M_-\to \bold 1$ be the counit. 
  
Define
the action of $A_+$ in $M_-$ by $\pi_-=m$; 
this is the same as the standard left action
of $A_+$ on itself. The coaction of $A_+$ in $M_-$ is then 
completely determined (via formula (2.3))
by its composition with $1_-$, which we
define by $\pi^*\circ 1_-=\iota\o 1_-$. 
It is easy to see from (2.3) that 
this coaction has the form
$$
\pi_-^*=(m^{op}\o 1)\circ (S^{-1}\o\sigma_{23})\circ \Delta_3.\tag 2.6
$$

The formulas for the action
and coaction of $A_+$ in $\hat M_+^*$ 
have the form
$$
\pi_+=m_3\circ \sigma_{23}\circ (1\o S^{-1}\o 1)\circ
 (\Delta^{op}\o 1),
\pi_+^*=\Delta^{op}.\tag 2.7
$$

For any object $X\in \Cal M_{A_+}^{A_+}$ consider the map
$\theta: \Hom_{\Cal M_{A_+}^{A_+}}( M_-\o X_0,\hat M_+^*\o X)\to 
\Hom_{\Cal N}(X_0,X)$ given by $\theta(f)=(\e\o 1_X)\circ f\circ 
(\iota \o 1_{X_0})$. 

\proclaim{Lemma 2.3} The map $\theta$ is an isomorphism.
\endproclaim

\demo{Proof} 
By Frobenius reciprocity, for any two dimodules $X,Y$ over 
$A_+$,
$\Hom_{\Cal M_{A_+}^{A_+}}(M_-\o X,\hat M_+^*\o Y)=
\Hom_{\Cal M^{A_+}}(X,\hat M_+^*\o Y)=
\Hom_{\Cal N}(X,Y)$. This implies the Lemma. 
\enddemo

Denote by $\psi$ the morphism $\theta^{-1}(1_X)$. 
Let $\eta=\psi\circ (\iota \o 1_{X_0}): X_0\to \hat M_+^*\o X$. 

\vskip .1in

2.4. {\it Dequantization of Hopf algebras.}

In the next three sections we will give the construction of the functor 
$\widehat{DQ}$ of dequantization, and show that this functor is inverse 
to the functor $\widehat{Q}$. 

Let us explain the idea of the construction of $\widehat{DQ}$.
We will start with a Hopf algebra $A_+$ in some symmetric tensor category 
$\Cal N$,
which satisfies a certain condition called quasisymmetry.
We will assign to it a family of Hopf algebras $A_+(t)$
depending algebraically on $t\in\Bbb Q$ (all realized on the same object of 
$\Cal N$), such that $A_+(1)=A_+$, and $A_+(0)$ is cocommutative. 
Then $A_+(0)$ has a natural co-Poisson Hopf structure, obtained as the 
quasiclassical limit of coproducts in $A_+(t)$ as $t\to 0$. 
Since this construction works in any symmetric tensor category, 
it is defined by universal acyclic formulas, 
and so it defines a functor $G$ from the category of 
quasisymmetric Hopf algebras 
in $\Cal N$ to the category of co-Poisson Hopf algebras in $\Cal N$. 

Now we can consider the case 
when $\Cal N$ is the category of topologically free $k[[h]]$-modules, 
and $A_+$ is a QUE algebra. In this case $A_+$ happens to be 
quasisymmetric. So we can define the Poisson-Hopf algebra $A_+(0)$
as above, and 
$\a_+$ to be the set of primitive elements of $A_+(0)$.
The object $\a_+$ has a natural Lie bialgebra structure. 
Thus we establish a functor $\widehat{DQ}$ from the category of
QUE algebras over $k[[h]]$ to the category of Lie bialgebras 
over $k[[h]]$, which will later be shown to be an inverse to $\widehat{Q}$. 

Let us now explain the main technical point -- how to construct
$A_+(t)$ from $A_+$. For this purpose we will consider 
the category $\Cal M_{A_+}^{A_+}$ of dimodules over $A_+$. 
This category is a braided tensor category, which has two
remarkable objects -- $M_-$ and $M_+^*$. 
In Section 1.4 we have explained how 
to recover the Hopf algebra $A_+$ using these objects. 

To construct $A_+(t)$, we choose a 1-parameter subsemigroup $a(t)$ 
in $\overline{GT}$, and twist the category $\Cal M^{A_+}_{A_+}$ 
by $a(t)$, as explained in Section 2.2. 
The objects $M_-,M_+^*$, which are crucial in the construction 
of Section 1.4, are also present in the twisted category 
$(\Cal M_{A_+}^{A_+})^{a(t)}$, and have the same properties. 
Therefore, one can apply the recovering procedure of Section 1.4
to this twisted category, which will yield the Hopf algebra $A_+(t)$, 
as desired. 

Now we describe the details of this construction. 

Let $A_+\in\Cal N$ be a Hopf algebra.

Let $\Cal M$ be the full tensor subcategory of the category
$\Cal M_{A_+}^{A_+}$ generated by $M_-$, $\hat M_+^*$. 
That is, objects of $\Cal M$ are arbitrary tensor products 
of copies of $\hat M_+^*$ and $M_-$, and morphisms are homomorphisms 
of dimodules.  

Let $\I$ be a tensor ideal in $\Cal M$ generated by 
the morphisms $R_{VW}-1, V,W\in \Cal M$ 
where $R$ is the R-matrix defined by
(2.5), and assume that $\Cal M$ is separated and complete
in the topology defined by $\Cal I$. In this case, we call
$A_+$ a quasisymmetric Hopf algebra. 

If $A_+$ is a quasisymmetric Hopf algebra, then the category $\Cal M$ is 
quasisymmetric. Pick $a\in \overline{GT}$.

Define a Hopf algebra $A_+^a$ as follows. 
As an object of $\Cal N$, $A_+^a$ coincides with $A_+$. 
The product is defined by formula (1.3), where $\Phi$ is the associativity 
morphism of the category $\Cal M^a$, 
and $1_+: \hat M_+^*\to\bold 1$ is the 
counit of $A_+$.

The coproduct is defined by the formula (1.4), where $\Delta_0$ 
is the coproduct in $A_+$, and $J: M_-\o M_-\to M_-\o M_-$ is 
given by the formula
$$ 
 J=(1_+\o 1\o 1_+\o 1)\circ \Phi_{1,3,24}^{-1}\circ \Phi_{3,2,4} 
\circ (RR^{op})^{(t-1)/2}(R^{op})^{-1}
\circ \Phi_{2,3,4}^{-1}\circ \Phi_{1,2,34}\circ (\eta\o\eta)\tag 2.8
$$
(see (1.5)). 

Using similar arguments to those given in Part II of \cite {EK},
one shows that 
these operations define a structure of a Hopf algebra on $A_+^a$. 

Let $QS(\Cal N)$ be the category of quasisymmetric Hopf algebras 
in $\Cal N$. We have constructed a functor $G_a:  
QS(\Cal N)\to QS(\Cal N)$ such that $G_b\circ G_a$ is isomorphic
to $G_{ab}$. Thus, we have an action of the semigroup 
$\widehat{GT}$ on the category $QS(\Cal N)$. 

Now fix $a_0\in GT_0$ be any point. 
Let
 $a(t)$ be the 1-parameter
semigroup in $\overline{GT}$ which is defined by Theorem 2.2.

Let $\Cal M(t)$, $t\in\Q$, be the tensor category $\Cal M^{a(t)}$.  
Define a family of Hopf algebras $A_+(t)=A_+^{a(t)}$, $t\in\Q$. 
It is clear that
the operations in $A_+(t)$ depend polynomially on $t$ modulo any
finite power of the ideal $\Cal I$. 

\proclaim{Proposition 2.4} The Hopf algebra $A_+(0)$ is cocommutative. 
\endproclaim

\demo{Proof} A direct computation.
\enddemo

\proclaim{Corollary 2.5} The Hopf algebra $A_+(0)$ has a natural structure
of a co-Poisson Hopf algebra, defined by $\delta=\lim_{t\to 0}
(\Delta-\Delta^{op})/t$. 
\endproclaim

We call the co-Poisson Hopf algebra $A_+(0)$ 
in $\Cal N$ the dequantization of $A_+$. 

\vskip .1in

2.5. {\it Dequantization of quantized universal enveloping algebras.}

Let $\Cal N$ be the category of topologically free $K$-modules,  
and $A_+\in HA_0(K)$ be a quantized universal enveloping algebra.

\proclaim{Proposition 2.6} 
$A_+$ is quasisymmetric.
\endproclaim
  
\demo{Proof} Easy
\enddemo

Therefore, applying the construction of Section 2.4, we obtain 
a co-Poisson Hopf algebra $A_+(0)$. It is clear that the assignment
$A_+\to A_+(0)$ is a functor from the category of QUE algebras
to the category of co-Poisson Hopf algebras, since it  
is given by a universal construction in the sense of Chapter 1. 

According to Drinfeld's Proposition 3.7 of \cite{Dr2}, any cocommutative 
Hopf QUE algebra $B$ is equal to $U(\b)$, where $\b$ is the Lie algebra
of primitive elements in $B$. Let $\a_+$ be the Lie algebra 
of primitive elements in $A_+(0)$. It has a natural structure of a Lie 
bialgebra, induced by the co-Poisson Hopf structure on $A_+(0)$. 

Thus, we have assigned to any QUE algebra $A_+$ a Lie bialgebra $\a_+$,
such that $A_+$ is isomorphic to $U(\a_+)$ mod $h$. This assignment
is clearly a functor $HA_0(K)\to LBA_0(K)$. We will call it the functor 
of dequantization, and denote
it by ${\widehat {DQ}_a}$, or, shortly, ${\widehat {DQ}}$. 

\vskip .1in

2.6. {\it Invertibility of the quantization functor.} 

Let $M_1(k)$ be the set of Lie associators over $k$.
For any $a=(\l,f)\in \overline{GT}(k)$, $\Phi\in M_1(k)$
define $a\Phi$ by formula (2.1):
$$
a\Phi=\Phi f(e^{h\Omega_{12}},\Phi^{-1}e^{h\Omega_{23}}\Phi). 
\tag 2.9
$$

Let $\Phi$ be the Lie associator over $k$ fixed in Section 1.3. 
Denote by $f_\Phi(X,Y)$ the 
expression of the form
$$
f_\Phi(X,Y)=e^{P_\Phi(\ln X,\ln Y)},\tag 2.10
$$
such that $\Phi=f_\Phi(e^{ht_{12}},e^{ht_{23}})$.
For any $\l\in k$, denote by $\Phi_\l$ the expression
$$
\Phi_\l=f_\Phi(X^\l,Y^\l).\tag 2.11
$$
 
The following result is due to Drinfeld, \cite{Dr4}.

\proclaim{Proposition 2.7} 
Let $\Psi$ be a Lie associator over $k$. There exists a unique
1-parameter semigroup $a_\Psi\subset \overline{GT}(k)$, 
such that $a_\Psi(\l)\in GT_\l(k)$ for all $\l\in k$, and 
$a_\Psi(\l)\Psi=\Phi_\l$. 
\endproclaim

\proclaim{Proposition 2.8} The functors $\widehat {Q_\Phi}$, 
${\widehat {DQ}_{a_\Phi}}$ 
are quasi-inverse
to each other. That is, the functors $\widehat {Q_\Phi}\circ 
{\widehat {DQ}_{a_\Phi}}, 
{\widehat {DQ}_{a_\Phi}}\circ \widehat {Q_\Phi}$ are
isomorphic to the identity. 
\endproclaim

\demo{Proof} Let $\a_+\in LBA_0(K)$ be a Lie bialgebra, $A_+=\widehat{Q}(\a_+)$.
For any $\l\in k$, let $\a_+(\l)$ be the Lie bialgebra obtained
from $\a_+$ by multiplying the cocommutator by $\l$. Let $A_+(\l)$ be
the Hopf algebra obtained from $A_+$ as described in Section 2.5, using
the semigroup $a=a_\Phi$. As $a_\Phi(\l)\Phi=\Phi^\l$,
we get that $A_+(\l)$ coincides with $\widehat{Q}(\a_+(\l))$. 
Tending $\l$ to $0$,
we obtain that the Lie bialgebra ${\widehat {DQ}}(A_+)$ 
is naturally isomorphic to 
$\a_+$. Thus, ${\widehat {DQ}_{a_\Phi}}\circ \widehat {Q_\Phi}$ is 
isomorphic to the identity.

Now let us prove that $\widehat{Q}_\Phi\circ \widehat{DQ}_{a_\Phi}$ is
isomorphic to the identity. 
Let $a_\Phi(0)=(0,g_\Phi)$. Let $A_+\in HA_0(K)$, 
$\a_+=\widehat{DQ}_{a_\Phi}(A_+)$, 
$\Cal M$ the category of $A_+$-dimodules defined in Section 2.4, 
$R$ the R-matrix defined by (2.5), and $T=R^{op}R$. Set
$(\Phi^0,\sigma R^0)=g_\Phi\cdot (1,\sigma R)$, where
$\sigma$ is the permutation. We have
$\Phi^0=g_\Phi(T_{12},T_{23})$, $R^0=RT^{-1/2}$. 

Define the Casimir element by $\Omega=\ln T$. 
Introduce commutativity and 
associativity morphisms in $\Cal M$ by
$$
R'=R^0e^{\Omega/2}, \Phi'=f_\Phi(\Phi^0
e^{\Omega_{12}}(\Phi^0)^{-1},e^{\Omega_{23}})\Phi^0.
\tag 2.12
$$
By the definition of $\widehat{Q}$, this tensor category  
is obtained in the process of quantization of the Lie bialgebra
$\a_+':=\widehat {DQ_{a_\Phi}}(A_+))$. Therefore,  
our claim follows from the equalities $R'=R$ and $\Phi'=1$. 

The first
equlity is obvious. It remains to prove the equality $\Phi'=\Phi$. 
Using the expression for $\Phi^0$, we see that
$$
\Phi'=f_\Phi(g_\Phi(T_{12},T_{23})^{-1}T_{12}
g_\Phi(T_{12},T_{23}),T_{23})g_\Phi(T_{12},T_{23}).\tag 2.13
$$
So, it remains to show that $f_\Phi\circ_1 g_\Phi=1$, where the operation 
$\circ_1$ is defined by (2.2). 

The set $\Cal F$ of exponentials of formal 
Lie series equipped with the operation $\circ_1$, 
is a group. By the definition
of $g_\Phi$, we have in this group $g_\Phi\circ f_\Phi=1$. 
This implies the desired equality $f_\Phi\circ_1 g_\Phi=1$.
Thus, the QUE algebra 
$\widehat{Q}(\widehat{DQ}(A_+))$ is canonically isomorphic to $A_+$. 
The Proposition  is proved.
\vskip .1in

Thus, we have shown that the functor $\widehat{Q}_\Phi$ 
is an equivalence of categories
and proved Theorem 2.1.

\end